\documentclass[11pt,a4paper,reqno]{amsart}

\usepackage{lineno,hyperref}
\usepackage{mathptmx}
\usepackage{amssymb}
\usepackage{amsmath}
\usepackage{amsthm}
\usepackage{latexsym}
\usepackage{amsfonts}
\usepackage{float}
\usepackage{hyperref}
\usepackage{placeins}
\usepackage{graphicx}
\usepackage[dvipsnames]{xcolor}
\usepackage{tikz}
\usetikzlibrary{arrows.meta,
	positioning,
	quotes}

\newtheorem{thm}{Theorem}[section]
\newtheorem{cor}[thm]{Corollary}
\newtheorem{lem}[thm]{Lemma}
\newtheorem{prop}[thm]{Proposition}

\theoremstyle{definition}
\newtheorem{defn}[thm]{Definition}
\newtheorem{rem}[thm]{Remark}
\newtheorem{ex}[thm]{Example}

\newtheorem{as}[thm]{Assumption}

\renewcommand{\emptyset}{\varnothing}
\newcommand{\rmref}[1]{{\rm\ref{#1}}}

\newcommand{\braces}[1]{{\rm (}#1{\rm )}}

\newcommand{\la}{\lambda}
\newcommand{\veps}{\varepsilon}
\newcommand{\vphi}{\varphi}
\newcommand{\ol}{\overline}
\newcommand{\wt}{\widetilde}

\newcommand{\R}{\ensuremath{\mathbb R}}    
\newcommand{\N}{\ensuremath{\mathbb N}}    

\newcommand{\calC}{\mathcal C}

\newcommand{\calM}{\mathcal M}

\newcommand{\calS}{\mathcal S}         
\newcommand{\calT}{\mathcal T}

\newcommand{\bmat}[4]{\begin{bmatrix}#1 & #2\\#3 & #4\end{bmatrix}}
\newcommand{\bvek}[2]{\begin{bmatrix}#1\\#2\end{bmatrix}}
\newcommand{\bmattt}[9]{\begin{bmatrix}#1 & #2 & #3\\#4 & #5 & #6\\#7 & #8 & #9\end{bmatrix}}
\newcommand{\bvekkk}[3]{\begin{bmatrix}#1\\#2\\#3\end{bmatrix}}

\newcommand{\sbvek}[2]{\left[\begin{smallmatrix}#1\\#2\end{smallmatrix}\right]}

\newcommand{\linspan}{\operatorname{span}}

\newcommand{\rank}{\operatorname{rank}}

\renewcommand{\d}{\mathrm{d}}

\newcommand{\dist}{\operatorname{dist}}

\newcommand{\XX}{\mathbb X}
\newcommand{\YY}{\mathbb Y}

\usepackage[normalem]{ulem}

\numberwithin{equation}{section}
\allowdisplaybreaks

\title{Optimal control of port-Hamiltonian systems: energy, entropy, and exergy$^*$}

\author[F.\ Philipp, M.\ Schaller, K.\ Worthmann, T.\ Faulwasser and B.\ Maschke]{Friedrich Philipp$^{1}$, Manuel Schaller$^{1}$, Karl Worthmann$^{1}$, Timm Faulwasser$^{2}$ and Bernhard Maschke$^{3}$}
	\thanks{}
		\thanks{$^{1}$Optimization-based Control Group, Institute for Mathematics, Technische Universit\"at Ilmemau, Germany
		{\tt\small \{friedrich.philipp, manuel.schaller, karl.worthmann\}@tu-ilmenau.de}.} %
	\thanks{$^{2}$Institute of Energy Systems, Energy Efficiency and Energy Economics, TU Dortmund University, Germany
		{\tt\small timm.faulwasser@ieee.org}}%
	\thanks{$^{3}$Univ Lyon, Universit{\'e} Claude Bernard Lyon 1, CNRS, LAGEPP UMR 
		5007, France {\tt\small bernhard.maschke@univ-lyon1.fr}\\
$^*$This research was initiated during a research stay of FP and MS at the LAGEPP UMR CNRS 5007 in Lyon partially funded
by the French Embassy in Germany by means of a PROCOPE mobility grant (M.\ Schaller). FP and MS thank the University of Lyon and the work group LAGEPP for the warm hospitality. 
FP was funded by the Carl Zeiss Foundation within the project \textit{DeepTurb--Deep Learning in and from Turbulence} and by the free state of Thuringia and the German Federal Ministry for Education and Research (BMBF) within the project \textit{THInKI--Th\"uringer Hochschulinitiative für KI im Studium}.
KW gratefully acknowledges funding by the German Research Foundation (DFG; project number 507037103). 
.}

\begin{document}
\begin{abstract}
We consider irreversible and coupled reversible-irreversible nonlinear port-Hamil\-tonian systems and the respective sets of thermodynamic equilibria. In particular, we are concerned with optimal state transitions and output stabilization on finite-time horizons. 
We analyze a class of optimal control problems, where the performance functional can be interpreted as a linear combination of energy supply, entropy generation, or exergy supply. 
Our results establish the integral turnpike property towards the set of thermodynamic equilibria 
providing a rigorous connection of optimal system trajectories to optimal steady states.
Throughout the paper, we illustrate our findings 
by means of two examples: a network of heat exchangers and a gas-piston system.

\end{abstract}

\maketitle

\section{Introduction}
\noindent Port-Hamiltonian systems provide a highly structured framework for energy-based modeling, analysis, and control of dynamical systems \cite{Geoplex09,Jeltsema_FTSC_14}. 
A particular feature of port-Hamilto\-nian systems is that 
solutions satisfy an energy balance and the supplied energy is given as a product of the input and the colocated 
output. This motivates and enables the formulation and in-depth analysis of optimal control problems, e.g., output stabilization or setpoint transitions, in which the supplied energy is minimized: 
Whereas the choice of this cost functional is physically meaningful, 
it leads to singular optimal control problems. 
However, for linear port-Hamiltonian systems, 
the port-Hamiltonian structure can be exploited to establish regularity of the optimality system~\cite{faulwasser2023hidden} (potentially after adding a rank-minimal regularization term) 
and to analyze the asymptotics of optimal solutions with respect to \ the conservative subspace~
\cite{Schaller2020a,Faulwasser2021}. These results (partly) have been extended to infinite-dimen\-sional linear systems~\cite{Philipp2021} and to nonlinear reversible port-Hamil\-tonian systems~\cite{karsai2023manifold,soledad2024conditions}. Further, for infinite-di\-men\-sional linear quadratic optimal control of reversible port-Hamiltonian systems, we refer to the recent works \cite{reis2024linear} in the context of system nodes and to \cite{hastir2024linear}, where networks of waves are studied. 
However, when the thermodynamic properties of the control systems have to be taken into account, other corresponding potentials than the free energy, such as the internal energy and the entropy, have to be considered. In turn, the dynamics have to reflect both---the energy conservation \textit{and} the entropy creation due to the irreversible phenomena.

The Hamiltonian formulation of controlled thermodynamic systems is an active field of research with various classes of systems ranging from dissipative Hamiltonian (or metriplectic) systems ~\cite{Oettinger_PhysRevE_06,Hoang_JPC_11,Hoang_JPC_2012} for the so-called GENERIC framework, irreversible port-Hamiltonian systems~\cite{Ramirez_ChemEngSci13,ramirez2022overview}, via port-Hamiltonian systems defined on contact manifolds \cite{Eberard_RMP07,Favache_IEEE_TAC_09,Favache_ChemEngSci10}, to symplectic manifolds \cite{Entropy_2018}.
%
Different control design methods have been developed for these systems: 
stabilization based on control of either thermodynamic potentials such as the availability function~\cite{Alonso96,Ydstie02,Ruszkowski05,Wang_LHMNLC15} or the entropy creation~\cite{GarciaSandoval_ChemEngSci_16,GarciaSandoval_JPC_17}, shaping of the entropy creation for irreversible port-Hamilto\-nian systems~\cite{Ramirez_Automatica16}, and structure-preserving stabilizing feedback of contact Hamiltonian systems \cite{Ramirez_SCL13,Ramirez_IEEE_TAC_17}.

%

In view of applications, due to the strong nonlinearity, only stationary optimal control problems were considered for high-dimensional thermodynamic systems (arising, e.g., from discretizations of partial differential equations), see, e.g., \cite{Johannessen_Energy_04_OptContrEntropy, DeKo2000}. 
In \cite[Section 13.5.2]{Wilh2015}, the authors discuss \emph{highways in state space} as particular sets, in which (or close to which) entropy optimal solutions reside. 
For non-stationary 
problems, such a behavior is coined the \emph{turnpike property} \cite{trelatzuazua,tudo:faulwasser22a,DaGrStWo2014}. 


In this work, we are concerned with optimal control of 
coupled reversible-irreversible port-Hamiltonian systems. To this end, we formulate the dynamic problem of energy, entropy, and exergy optimization. 
We show that, for long time horizons, optimal solution trajectories of this problem stay close to the set of thermodynamic equilibria for the majority of the time, for both output stabilization and set point transitions as the underlying control task. In this context, we significantly extend our previous conference paper~\cite{MascPhil22} wherein only set point transitions of an irreversible system without reversible components is studied.
Further, we provide a detailed analysis of the underlying steady-state optimization problem
and we illustrate the theoretical results by means of several numerical examples including a network of heat exchangers and a gas-piston system.



This paper is structured as follows: In Section~\ref{sec:setting}, we \hbox{recall} the definition of the class of reversible-irreversible port-Ha\-mil\-to\-nian systems 
and show that two examples embed into this framework: 
a heat exchanger and a gas-piston system. We provide conditions in terms of the Hamiltonian energy, which ensure that the set of thermodynamic equilibria is a manifold and characterize its dimension.
Next, in Section~\ref{sec:statetrans}, we consider and analyze the problem of optimal state transition, where optimality is understood as minimal energy supply, minimal entropy creation, or a combination of both as the exergy. In Section~\ref{sec:outputstab}, we derive similar results for optimal output stabilization instead of state transition, demonstrating the generality of the developed tools. Finally, we illustrate our findings by reconsidering the two examples in numerical illustrations in Section~\ref{sec:numerics} before conclusions are drawn in Section~\ref{sec:conclusions}.



\medskip
\noindent \textbf{Notation}. We denote the gradient of a scalar valued function $F:\R^n \to \R$ by $F_x$ and its Hessian by $F_{xx}$. The {\em Poisson bracket} with respect to a skew-symmetric matrix $J\in\R^{n\times n}$ of the differentiable functions $v,w : \R^n\to\R$ is defined by
$$
\{v,w\}_J(x) := v_x(x)^\top Jw_x(x),\quad x\in\R^n.
$$
Let $K\subset\R^n$. We write ``$\alpha(x)\lesssim\beta(x)$ for $x\in K$'' if there exists $c>0$ such that $\alpha(x)\le c\beta(x)$ for $x\in K$. In addition, we write ``$\alpha(x)\asymp\beta(x)$ for $x\in K$'' if $\alpha(x)\lesssim\beta(x)$ and $\beta(x)\lesssim\alpha(x)$ for $x\in K$.

\section{Coupled reversible-irreversible pH-Systems}\label{sec:setting}\noindent In this section we provide the system class we analyze and illustrate it by means of two examples.

\subsection{Definition of reversible-irreversible pH-Systems}
\noindent The state space will be given by an open set $\mathbb X\subset\R^n$, $n\in\N$. As usual in finite-dimensional pH systems, input and output space coincide; here, inputs and outputs are elements of $\R^m$, $m\in\N$. The following definition of coupled reversible-irreversible pH Systems was given in \cite{Ramirez_EJC13}.

\begin{defn}[Reversible-irreversible pH-Systems]\label{def:RIPHS}~\\
A \emph{Reversible-Irreversible port-Hamiltonian system} (RIPHS) is defined by 
\begin{itemize}
\item[(i)] a locally Lipschitz-continuous map $J_0 : \XX\to\mathbb{R}^{n\times n}$ of skew-symmetric \emph{Poisson structure matrices} $J_0(x)$, $x\in\XX$,
\item[(ii)]  a differentiable \emph{Hamiltonian  function} $H : \XX\to\R$, whose gradient $H_x : \XX\to\R^n$ is locally Lipschitz-continuous,
\item[(iii)] a differentiable \emph{entropy function} $S :  \XX\to\R$, whose gradient $S_x : \XX\to\R^n$ is locally Lipschitz-continuous and which is a Casimir function of the Poisson structure matrix function $J_0$, that is, $J_0(x)S_x(x)=0$ for all $x\in\XX$,
\item[(iv)] $N$ constant \emph{skew-symmetric structure matrices} $J_{k}\in\mathbb{R}^{n\times n}$, $k=1,\ldots, N$,
\item[(v)]  functions $\gamma_1,\ldots,\gamma_N : \R^{2n}\to\R$ such that $\gamma_{k}\left(\,\cdot\,,H_x(\cdot)\right) : \XX\to\R$ is positive and locally Lipschitz-continuous for all $k=1,\ldots, N$,
\item[(vi)] a locally Lipschitz-continuous input matrix function $g : \R^{2n}\to\R^{n\times m}$,
\end{itemize}
the \emph{state equation}
\begin{equation}\label{eq:RIPHS}\tag{RIPHS}
\dot{x}=\Big(J_{0}(x)+\sum_{k=1}^{N}\gamma_{k}(x,H_x)\left\{ S,H\right\} _{J_{k}}J_{k}\Big)H_x(x)+g(x,H_x)u,
\end{equation}
and two conjugate outputs, corresponding to the energy and the entropy, respectively, defined by
\begin{align}\label{eq:outputs}
y_H := g(x,H_x)^\top H_x \qquad\text{and}\qquad y_S := g(x,H_x)^\top S_x.
\end{align}
\end{defn}
\noindent If $J_0\equiv 0$ and $N=1$, the system is called an {\em irreversible port-Hamiltonian system} (IPHS).

By means of direct calculations, any $x\in C^1(0,T;\XX)$ satisfying the dynamics \eqref{eq:RIPHS} can be shown to obey both the {\em energy balance}
\begin{align}\label{eq:energybalance_sys}
\frac{\text{d}}{\text{d}t}H(x(t)) &= y_H(t)^\top u(t)
\end{align}
and the {\em entropy balance}
\begin{align}\label{eq:entropybalance}
\frac{\text{d}}{\text{d}t}S(x(t)) &= \sum_{k=1}^N \gamma_k\big(x(t),H_x(x(t))\big) \{S,H\}_{J_k}^2(x(t)) + y_S(t)^\top u(t)
\end{align}
Here, the quantity $y_H(t)^\top u(t)$ in \eqref{eq:energybalance_sys} represents the {\em energy flow}, i.e., the power supplied to the system, whereas $y_S(t)^\top u(t)$ in \eqref{eq:entropybalance} stands for the {\em entropy flow} injected into the system. The positivity of the sum on the right-hand side of the entropy balance captures the irreversible nature of the particular phenomenon. For closed systems, i.e., where $g(x,H_x)\equiv 0$, it can be directly read off these equations that energy is conserved and entropy is non-decreasing, hence the two fundamental laws of thermodynamics are satisfied.

For $u=0$, a \emph{thermodynamic equilibrium} is attained if the first term on the right-hand side of \eqref{eq:entropybalance} vanishes. The set of {\em thermodynamic equilibria}  
\begin{align*}
\mathcal{T} :=& \Big\{x\in\XX : \sum_{k=1}^{N} \gamma_k(x,H_x(x))\{S,H\}^2_{J_k}(x)=0\Big\}
\\=& \big\{x\in\XX : \{S,H\}_{J_k}(x)=0\,\, \text{for all }k=1,\ldots,N\big\},
\end{align*}
plays a distinguished role in this paper. In the set definition above, the second equality follows from positivity of $\gamma$. Out of thermodynamic equilibrium, the energy balance equation \eqref{eq:energybalance_sys} reduces to the conservation of energy and the entropy balance equation \eqref{eq:entropybalance} reduces to the irreversible entropy creation.

\begin{rem}[Linear and affine input structures]
Note that the input map in Definition~\ref{def:RIPHS} is \emph{linear} in the input, contrary to the setting in \cite{Ramirez_EJC13} where it is affine. This choice  corresponds to a reversible interconnection of the system with its environment as discussed in \cite{Maschke_submIFAC_WC23b}.
\end{rem}

\begin{rem}[Local existence of solutions]
Given a control $u\in L^\infty_{\rm loc}((t_0,t_f),\R^m)$, the local Lipschitz conditions on the functions $J_0$, $H_x$, $S_x$, $\gamma_k$, and $g$ guarantee at least local existence and uniqueness of solutions of initial value problems with dynamics \eqref{eq:RIPHS}.
\end{rem}

In the next lemma, we assume that the entropy $S$ is linear in the state $x$. This is not a restrictive assumption as in many physical examples, where the Hamiltonian~$H$ contains the internal energies of subsystems, the total entropy is simply the sum of some coordinates of the state (see Example \ref{ex:gp} below).


\begin{lem}\label{l:manifold}
The set of thermodynamic equilibria $\calT$ is closed in $\XX$. If $S(x) = e^\top x$ with some $e\in\R^n$ and for all $x\in\XX$, we have
\begin{equation}\label{e:suff}
\rank\sbvek{H_{xx}(x)}{H_x(x)^\top} = n,
\end{equation}
then $\calT$ is a $C^1$-submanifold of $\XX$. It is empty if and only if $H_x(\XX)\cap\bigcap_{k=1}^N(J_ke)^\perp = \emptyset$. Otherwise, its dimension is given by
$$
\dim\calT = n - \rank[J_1e\cdots J_Ne].
$$
\end{lem}
\begin{proof}
As the Poisson bracket $\{S,H\}_{J_k}$ is continuous for any $k\in\{1,\ldots,N\}$, it is clear that $\calT$ is closed in $\XX$. Now, assume that $S(x) = e^\top x$ and that \eqref{e:suff} holds. Set $v_k = J_ke$, $k=1,\ldots,N$. Note that
$$
\{S,H\}_{J_k} = S_x^\top J_kH_x = e^\top J_kH_x = -v_k^\top H_x.
$$
If $v_1=\ldots=v_N=0$, then $\calT = H_x^{-1}(\R^n) = \XX$ and thus $\dim\calT = n$, as claimed. Otherwise, $r:=\rank[v_1\dots v_n]\ge 1$. We assume without loss of generality that $v_1,\ldots,v_r$ are linearly independent and define $f : \XX\to\R^r$ by
$$
f(x) = [v_1,\dots,v_r]^\top H_x(x).
$$
Note that $\calT = f^{-1}(\{0\})$. We have $Df(x) = [v_1,\dots,v_r]^\top H_{xx}(x)$, $x\in\XX$. Let $x_0\in\calT$, i.e., $H_x(x_0)^\top v_k=0$ for $k=1,\ldots,r$ and suppose that $Df(x_0)$ is not surjective, i.e., there exists some $w\in\R^r\backslash\{0\}$ such that $H_{xx}(x_0)[v_1,\dots v_r]w=0$. Consequently,
$$
\sbvek{H_{xx}(x)}{H_x(x)^\top}[v_1,\dots v_r]w = 0,
$$
so that, by assumption, $[v_1,\dots v_r]w = 0$. But $v_1,\dots,v_r$ are linearly independent so that $w=0$, a contradiction. Therefore, zero is a regular value of the $C^1$-function $f$, hence $\calT = f^{-1}(\{0\})$ is a $C^1$-submanifold of $\XX$ of the given dimension, cf.\ \cite[Theorem 73C, p.\ 556]{Zeidler88}.
\end{proof}

We close this introductory subsection with the observation that the class of reversible-irreversible port-Hamiltonian systems is invariant under linear coordinate transforms.

\begin{prop}
If $x$ solves \eqref{eq:RIPHS} and $V\in\R^{n\times n}$ is invertible, then $z=Vx$ solves the reversible-irreversible pH-system
\begin{align*}
\dot{z}=\Big(\wt J_{0}(z)+\sum_{k=1}^{N}\wt\gamma_k(z,\wt H_z)\big\{\wt S,\wt H\big\}_{\wt J_{k}}\wt J_{k}\Big)\wt H_z(z) + \tilde g(z,\wt H_z)u
\end{align*}
with energy Hamiltonian $\wt H(z) = H(x)$, entropy $\wt S(z) = S(x)$,
$$
\wt J_0(z) = VJ_0(x)V^\top,\quad \wt J_k = VJ_kV^\top,
$$
and
$$
\wt g(z,\wt H_z) = Vg(x,V^\top\wt H_z),\quad
\wt\gamma_k(z,\wt H_z) = \gamma_k(x,V^\top\wt H_z).
$$
\end{prop}
\begin{proof}
The main observation is that, setting, $\wt H_z(z) = V^{-\top}H_x(x)$ and $\wt S_z(z) = V^{-\top}S_x(x)$, we have
\begin{align*}
\gamma_{k}(x,H_x)\left\{ S,H\right\}_{J_{k}}
&= \gamma_{k}(x,V^\top\wt H_z)S_x^\top J_k H_x = \wt\gamma_{k}(z,\wt H_z)\{\wt S,\wt H\}_{\wt J_k}.
\end{align*}
Thus
\begin{align*}
\dot z
&= V\Big(J_{0}(x)+\sum_{k=1}^{N}\wt\gamma_{k}(z,\wt H_z)\{\wt S,\wt H\}_{\wt J_k}J_{k}\Big)V^\top \wt H_z(z)+ Vg(x,H_x)u\\
&= \Big(\wt J_{0}(z)+\sum_{k=1}^{N}\wt\gamma_{k}(z,\wt H_z)\{\wt S,\wt H\}_{\wt J_k}\wt J_{k}\Big)\wt H_z(z)+ \wt g(z,\wt H_z)u,
\end{align*}
which proves the proposition.
\end{proof}

\subsection{Tutorial Examples}
\noindent Next, we illustrate the class of \eqref{eq:RIPHS} systems with two examples from thermodynamics:  a heat exchanger and a gas-piston system.

\begin{ex}[Heat exchanger]\label{ex:heat}
Let us consider a very simple model of a heat exchanger as depicted in Figure~\ref{fig:heat}. The example is slightly adapted from \cite{RamirezThesis:2012}. The variables $T_i$,$ S_i$, and $H_i$ denote the temperature, the entropy and the energy in compartment $i=1,2$, respectively.  
\begin{figure}[ht]
\includegraphics[width=.6\linewidth]{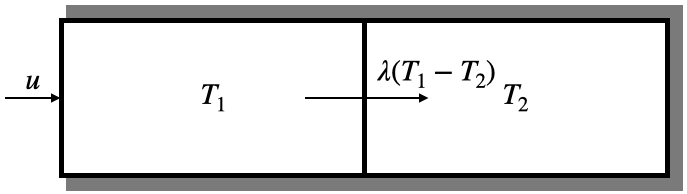}
\caption{Simple model of a heat exchanger}\label{fig:heat}
\end{figure}
Assuming that the walls are non-deformable and impermeable, the thermodynamic properties of each compartment are given by the following relation between the temperature and the entropy
\begin{align}\label{eq:relation}
    T_i(S_i) = T_\text{ref}\cdot e^{(S_i-S_\text{ref})/c_i}, \qquad i=1,2,
\end{align}
where $S_\text{ref}\in \mathbb{R}$ is a reference entropy corresponding to the reference temperature $T_\text{ref}$ and $c_i$, $i=1,2$, are heat capacities, cf.\ \cite[Section 2.2]{Couenne06}. The energy in each compartment, denoted by $H_i(S_i)$, $i=1,2$, can be obtained by integrating Gibbs' equation $\text{d}H_i = T_i\text{d} S_i$ as a primitive of the function $T_i(S_i)$, $i=1,2$. The state of our IPHS is $x=\begin{bmatrix}S_1&S_2\end{bmatrix}^\top$ and the total energy (entropy) is given by the sum of the energies (entropies, resp.) in the compartments, i.e.,
\begin{align*}
H(x) = H_1(S_1) + H_2(S_2) \quad\text{and}\quad S(x) = S_1+S_2 = \begin{bmatrix}1&1\end{bmatrix}x.
\end{align*}
We first assume that the outer walls are perfectly isolated, i.e., there is only heat flux through the conducting wall in between the two compartments given by Fourier's law
\begin{align}
\label{e:fourier}
q = \lambda(T_1-T_2).
\end{align}
On the other hand, the change in energy is given by the heat flux and thus, by continuity of the latter, we have
\begin{align*}
q = -\tfrac{\text{d}}{\text{d}t} H_1(S_1(t)) = \tfrac{\text{d}}{\text{d}t}H_2(S_2(t)).
\end{align*}
By Gibb's equation we have $\frac{\d}{\d t}(H_j\circ S_j) = \frac{\d H_j}{\d S_j}\cdot\frac{\d}{\d t}S_j = T_j\cdot\frac{\d}{\d t}S_j$. Hence,
\begin{align*}
\lambda(T_1-T_2)=-T_1\tfrac{\text{d}}{\text{d}t}S_1(t) = T_2\tfrac{\text{d}}{\text{d}t}S_2(t).
\end{align*}
Rearranging this equation yields the Hamiltonian-like formulation
\begin{align*}
\tfrac{\text{d}}{\text{d}t} \bvek{S_1(t)}{S_2(t)}
= \lambda\big(\tfrac{1}{T_2(t)}-\tfrac{1}{T_1(t)}\big)\underbrace{\bmat 0{-1}1{\phantom{-}0}}_{=:J}\bvek{T_1(t)}{T_2(t)}.
\end{align*}
To complete the definition of an uncontrolled IPHS, we set
\begin{align*}
\gamma(x,H_x) = \tfrac{\lambda}{T_1T_2}.
\end{align*}
Then, as
$$
\qquad \{S,H\}_J = \big[1\quad 1\big]\,J \bvek{T_1}{T_2} = T_1-T_2,
$$
the above ODE is of the form \eqref{eq:RIPHS} with $J_0\equiv 0$, $N=1$, and $g(x,H_x)\equiv 0$.

\smallskip
\noindent\textbf{Entropy flow control.} The first possible choice of an input would be to control the entropy flow $u$  into or out of compartment one, cf.\ Figure~\ref{fig:heat}, cf.\ \cite[Section 4.4.4]{RamirezThesis:2012}. In this case, we have $g(x,H_x)\equiv\sbvek 10$ and hence the dynamics
\begin{align}\label{eq:entcont}
\frac{\text{d}}{\text{d}t} \bvek{S_1}{S_2}  &= \gamma(x,H_x) \{S,H\}_J(x) J H_x + \bvek 10 u.
\end{align}

\noindent\textbf{Control by thermostat.} A choice which is realizable in practise is to connect the first compartment to a thermostat with a controlled temperature $T_e$. Consequently, the heat flow between this compartment and the thermostat can be described via
\begin{align*}
q_e = \lambda_e(T_e-T_1),
\end{align*}
with $\lambda_e>0$ being a heat conduction coefficient.

In this case, the entropy flow control \eqref{eq:entcont} is related to the thermostat temperature control $T_e$  by the state dependent control transformation
\begin{equation}\label{eq:NonlinearControlTransformation}
 u = \tfrac{T_e - T_1}{T_1}   
\end{equation}

\end{ex}

\begin{ex}[Gas-piston system, cf.\ \cite{vdS23,Ramirez_EJC13}]\label{ex:gp}
Let us briefly recall the model of a gas contained in a cylinder closed by a piston subject to gravity. Contrary to the previous example of the heat exchanger which was purely thermodynamic, the gas-piston system will encompass the thermodynamic and the mechanical domain.

\smallskip
\noindent\textbf{Energy and co-energy variables}.
Consider the gas in the cy\-lin\-der under the piston and assume that it is closed, i.e., there is no leakage.  Then, the internal energy $U$ of the gas may be expressed as a function of its entropy $S$ and its volume $V$. For an ideal gas (see \cite{Eberard_RMP07}), we have
\begin{equation}\label{eq:InternalEnergyIdealGas}
U(S,V)=\frac{3}{2}\,NRT_{0}\cdot e^{\beta(S,V)},
\end{equation}
with the ideal gas constant $R$ and 
\begin{equation}
\beta(S,V)=\frac{S-Ns_{0}+RN{\rm ln}(NRT_{0})-RN{\rm ln}(VP_{0})}{\frac{3}{2}RN}\label{eq:InternalEnergyIdealGasAdd}
\end{equation}
where $s_{0},\,T_{0},P_{0}$ are positive reference values.
Assume that the energy of the piston consists of the sum of the gravity potential energy and kinetic energy:
$$
H_{mec} = \frac{1}{2m}\,p^{2}+mgz,
$$
where $z$ denotes the altitude of the piston, $p$ its kinetic
momentum, and $m$ its mass. If we assume that the gas fills all the volume below of the piston, then we have $Az = V$ where $A$ stands for the base area of the piston's head. Hence one may choose the vector of  state or energy variables $x=[S,V,p]^{\top}$ and the total energy of the system is given by
$$
H(x)=U(S,V)+H_{mec}(V,p)
$$
where
$$
H_{mec}(V,p) = \frac{1}{2m}\,p^{2}+\frac{mg}{A}V.
$$
The differential of the energy defines the co-energy variables
\begin{subequations}\label{eq:Co-energyVariablesPistonGas}
\begin{align}
\frac{\partial H}{\partial S}&= \frac{\partial U}{\partial S}=T_{0}\,e^{\beta(S,V)} = T(S,V)\\
\frac{\partial H}{\partial V}&=\frac{\partial U}{\partial V} + \frac{\partial H_{\rm mec}}{\partial V}  = -P(S,V)+ \frac{mg}{A} = \frac{-NRT_{0}\,e^{\beta(S,V)}}{V} + \frac{mg}{A} \\
\frac{\partial H}{\partial p}&=\frac{\partial H_\text{mec}}{\partial p} = \frac{p}{m} = v(p)
\end{align}
\end{subequations}
where $T$ is the temperature, $P$ is the pressure, 
and $v$ is the velocity of the piston.

\smallskip
\noindent\textbf{Dynamics: reversible and irreversible coupling.}
The model consists of three coupled balance equations. It may be written as a \textit{quasi-Hamiltonian} system with a skew-sym\-metric structure matrix depending on two co-energy variables, the velocity $v$ and the temperature $T$:
\begin{equation}
\frac{d}{dt}
\bvekkk SVp
=
\underbrace{\bmattt {\phantom{-}0}{\phantom{-}0}{\frac{\kappa v}{T}}{\phantom{-}0}{\phantom{-}0}A{-\frac{\kappa v}{T}}{-A}0}_{=:J_{irr}\left(T,v\right)}
\underbrace{\bvekkk T{-P+ \frac{mg}{A}}v}_{H_x}\label{dynamic_gas-piston}.
\end{equation}
The first equation is the entropy balance accounting for the irreversible creation of entropy associated with the \textit{non-reversible phenomena} due to mechanical friction and viscosity of the gas when the piston moves. The second equation relates the variation of the volume of the gas to the velocity of the piston. The last equation is the momentum balance equation accounting for the mechanical forces induced by gravity, the pressure of the gas, and the total resistive forces which are assumed to be linear in the velocity of the piston, i.e., $F_{r}(v)=\kappa v$ with a scalar constant $\kappa\geq 0$.

The system \eqref{dynamic_gas-piston} may be written in the form of \eqref{eq:RIPHS} by decomposing further the structure matrix
\[
J_{irr}\left(T,v\right)=J_{0}+R\left(x,H_x,S_x\right)J_1
\]
with the constant \emph{Poisson} structure matrix associated with the reversible coupling
\[
J_0 = \bmattt 0{\phantom{-}0}00{\phantom{-}0}A0{-A}0
\]
Since $S_x(x) = [1,0,0]^\top$, we have $J_0S_x = 0$, i.e., the entropy is a Casimir function of $J_0$, as required in Definition~\ref{def:RIPHS}(iii). The second structure matrix corresponding to the dissipative phenomenon, that is, the friction of the piston is given by
\[
J_1=\bmattt {\phantom{-}0}01{\phantom{-}0}00{-1}00.
\]
The modulating function for the irreversible phenomenon is
$$
R(x,H_x,S_x) = \gamma(x,H_x)\{S,H\}_{J_1},
$$
which is composed by the Poisson bracket 
\[
\{S,H\}_{J_1} = \begin{bmatrix}1 & 0 & 0\end{bmatrix}J_1\bvekkk {\phantom{-}T}{-P}{\phantom{-}v} = v,
\]
which is indeed the driving force inducing the mechanical friction and viscosity forces in the gas, and the positive function
\[
\gamma\left(x,H_x\right)=\gamma\left(T\right)=\frac{\kappa}{T}.
\]


\noindent \textbf{Control input.}
For the purpose of illustration, we assume that the system is controlled by a heat flow $u$ entering the first line of the dynamics \eqref{dynamic_gas-piston}, similar to the entropy flow in \eqref{eq:entcont}. Similar to the heat exchanger it would be more realistic to assume that this heat flow is generated by a thermostat at controlled temperature $T_{e}$, interacting with the gas through a heat conducting wall leading to the same control transformation~\eqref{eq:NonlinearControlTransformation}.

\end{ex}

\section{Optimal state transitions}\label{sec:statetrans}
\noindent Let a horizon $t_f\geq 0$, an initial state $x^0\in\XX$ and a terminal region $\Phi \subset\XX$ be given and assume that $\mathbb{U}\subset\mathbb{R}^m$ is compact and convex, containing the origin in its interior, i.e., $0\in\operatorname{int}\mathbb U$. Let $\alpha_1,\alpha_2 > 0$ and consider the prototypical optimal control problem  
\begin{align}
\label{eq:phOCP}\tag{phOCP}
\begin{split}
\min_{u\in L^1(0,t_f;\mathbb{U})} \int_0^{t_f}&\left[\alpha_1y_H(t) - \alpha_2 T_0y_S(t)\right]^\top u(t)\,\text{d}t \\
\text{s.t. } \eqref{eq:RIPHS},\quad &x(0)=x^0, \quad x(t_f)\in \Phi ,\quad \eqref{eq:outputs},
\end{split}
\end{align}
where $T_0>0$ is a reference temperature and $L^1(0,t_f;\mathbb U)$ denotes the space of all measurable functions $u : [0,t_f]\to\mathbb U$. Note that $L^1(0,t_f;\mathbb{U}) \subset L^\infty(0,t_f;\mathbb{U})$ due to boundedness of $\mathbb{U}$. From this optimal control problem, we can recover three important choices for the cost functional:
\begin{itemize}
\item $\alpha_1 = 1$, $\alpha_2=0$: Minimal energy supply.
\item $\alpha_1 = 0$, $\alpha_2=1$: Minimal entropy extraction.
\item$\alpha_1 = 1$, $\alpha_2=1$: Minimal exergy supply.
\end{itemize}
Setting $\ell_{\alpha_1,\alpha_2}(x,u) = [\alpha_1y_H - \alpha_2 T_0 y_S]^\top u$ and using the balance equations \eqref{eq:energybalance_sys} and \eqref{eq:entropybalance}, we obtain the identity
\begin{align}
\begin{split}\label{eq:energybalance}
\int_0^{t_f} \ell_{\alpha_1,\alpha_2}(x,u)\,\text{d}t &= \alpha_1\left[H(x(t_f))-H(x^0)\right]\\ &\qquad + \alpha_2T_0\left(\!S(x^0)- S(x(t_f)) +\! \sum_{k=1}^N\int_0^{t_f}\!\!\!\gamma_k(x,H_x)\{S,H\}^2_{J_k}\,\text{d}t\right).
\end{split}
\end{align}



\noindent 
%

\subsection{Optimal steady states}
\noindent In this part, we analyze the steady-state counterpart of the optimal control problem \eqref{eq:phOCP}. In the literature considering entropy optimization, in particular in the context of distributed parameter systems, the steady state problem is often considered due to the high complexity of the dynamic problem, cf.\ e.g.\ \cite{Johannessen_Energy_04_OptContrEntropy} for a plug flow reactor model and~\cite{DeKo2000} for a binary tray distillation process.

In the turnpike result established in Subsection\ref{subsec:statetrans_turnpike}, we provide a rigorous connection between dynamic and static problem, stating that entropy-, energy-, and entropy-optimal solutions for the dynamic problem are close to optimal solutions of the steady state for the majority of the time. This in particular consolidates the central role of the steady state problem in the context of (dynamic) optimal control for irreversible systems.

%
\noindent To simplify notation, we write
$$
f(x,u) = \wt J(x)H_x(x) + g(x,H_x)u,
$$
where
\begin{align}\label{e:tildeJ}
\wt J(x) = J_0(x) + \sum_{k=1}^N\gamma_k(x,H_x)\cdot\{S,H\}_{J_k}\!(x)\,J_k.
\end{align}
Then \eqref{eq:RIPHS} reads $\dot x = f(x,u)$. A {\em steady state} of \eqref{eq:RIPHS} is a pair $(\bar x,\bar u)\in\XX\times\mathbb U$ for which $f(\bar x,\bar u) = 0$. In particular, the constant trajectory $x(t)\equiv\bar x$ is a solution to $\dot x = f(x,\bar u)$.

The steady state optimal control problem corresponding to \eqref{eq:phOCP} reads as follows:
\begin{align}\label{e:OSSP}
\min_{(\bar x,\bar u)\in\XX\times\mathbb U}\ell_{\alpha_1,\alpha_2}(\bar x,\bar u)\qquad\text{s.t.}\quad f(\bar x,\bar u)=0.
\end{align}
Any solution of this problem will be called an optimal steady state for \eqref{eq:phOCP}.

\begin{prop}\label{p:oss}
If $(\bar x,\bar u)$ is a steady state of \eqref{eq:RIPHS}, then
\begin{equation}\label{e:ge0}
\ell_{\alpha_1,\alpha_2}(\bar x,\bar u) = \alpha_2 T_0\cdot\sum_{k=1}^N\gamma_k(\bar x,H_x(\bar x))\cdot\{S,H\}_{J_k}^2(\bar x)\ge 0.
\end{equation}
In particular, the optimal value of \eqref{e:OSSP} is always non-negative and it is zero if and only if the set
$$
\calS = \{(\bar x,\bar u)\in\calT\times\mathbb U : g(\bar x,H_x(\bar x))\bar u = -J_0(\bar x)H_x(\bar x)\}
$$
is non-empty. In this case, the set $\calS$ coincides with the set of optimal steady states for \eqref{eq:phOCP}. If $H_x(\hat x)\in\linspan\{S_x(\hat x)\}$ for some $\hat x\in\XX$, we have $(\hat x,0)\in\calS$.
\end{prop}
\begin{proof}
Assuming $f(\bar x,\bar u) = 0$, we compute
\begin{align*}
\ell_{\alpha_1,\alpha_2}(\bar x,\bar u)
&= \alpha_1\cdot y_H^\top \bar u - \alpha_2T_0\cdot y_S^\top \bar u\\
&= \big[\alpha_1\cdot H_x(\bar x) - \alpha_2T_0\cdot S_x(\bar x)\big]^\top g(\bar x,H_x(\bar x))\bar u\\
&= -\big[\alpha_1\cdot H_x(\bar x) - \alpha_2T_0\cdot S_x(\bar x)\big]^\top\wt J(\bar x)H_x(\bar x)\\
&= \alpha_2T_0\cdot S_x(\bar x)^\top\wt J(\bar x)H_x(\bar x).
\end{align*}
Now, we substitute $\wt J(\bar x)$ with the right hand side from \eqref{e:tildeJ}, further we note that $S_x(\bar x)^\top J_0(\bar x) = 0$ (see Definition \ref{def:RIPHS} (iii)) and obtain \eqref{e:ge0}.

Let $(\hat x,\hat u)\in\calS$. Then $J_0(\hat x)H_x(\hat x) + g(\hat x,H_x(\hat x))\hat u = 0$, and, as $\hat x\in\calT$, we also have $\{S,H\}_{J_k}(\hat x) = 0$ for all $k=1\ldots,N$. Hence $f(\hat x,\hat u) = 0$, which means that $(\hat x,\hat u)$ is a steady state of \eqref{eq:RIPHS}. The first part of the proposition hence further implies that $\ell_{\alpha_1,\alpha_2}(\hat x,\hat u) = 0$, which shows that the optimal value of \eqref{e:OSSP} is indeed zero and that $(\hat x,\hat u)$ is an optimal steady state for \eqref{eq:phOCP}. If $(\bar x,\bar u)$ is any other optimal steady state for \eqref{eq:phOCP}, then \eqref{e:ge0} implies that also $\{S,H\}_{J_k}\!(\bar x) = 0$ for $k=1,\ldots,N$. That is, $(\bar x,\bar u)\in\calS$.

If the optimal value of \eqref{e:OSSP} is zero and $(\bar x,\bar u)$ is a minimizer, then \eqref{e:ge0} implies that $\bar x\in\calT$ and thus $(\bar x,\bar u)\in\calS$.

Now suppose there exists some $\hat x\in\R^n$ such that $H_x(\hat x)\in\linspan\{S_x(\hat x)\}$. Let $H_x(\hat x) = \kappa\cdot S_x(\hat x)$ for some $\kappa\in\R$. Then we have $\{S,H\}_{J_k}(\hat x) = S_x(\hat x)^\top J_k H_x(\hat x) = 0$ for all $k=1,\ldots,N$, thus $\hat x\in\calT$. Moreover, $J_0(\hat x)H_x(\hat x) = 0$ by Definition \ref{def:RIPHS} (iii). This implies $(\hat x,0)\in\calS$.
\end{proof}

\begin{cor}\label{c:irr}
In the irreversible case, i.e., $J_0\equiv 0$ and $N=1$, we have $\calT\times\{0\}\subset\calS$. In particular, if $\calT$ is non-empty, then so is $\calS$ and thus $\calS$ coincides with the set of all optimal steady states for \eqref{eq:phOCP}. 
\end{cor}

\noindent Having defined the sets $\mathcal{T}$ of thermodynamic equilibria, and the set $\mathcal{S}$ being closely related to optimal steady states, a third set comes into play, which is the set of thermodynamic equilibria which can be losslessly maintained while obeying the dynamics of \eqref{eq:RIPHS}:
$$
\calT_{\rm opt} := \{\bar x\in\calT : \exists\bar u\in\mathbb U \text{ with } (\bar x,\bar u)\in\calS\}.
$$
%
Note that in the irreversible case we have $\calT_{\rm opt} = \calT$.

\subsection{Turnpikes towards thermodynamic equilibria}
\label{subsec:statetrans_turnpike}
\noindent In this subsection we shall impose the following assumptions on the Hamiltonian $H$ and the entropy function $S$.

\begin{as}\label{as:ham}
\phantom{foo}
\begin{enumerate}
\item[(a)] $H\in C^2(\XX,\R)$, the image $\YY := H_x(\XX)$ is open in $\R^n$, and $H_x: \XX\to\YY$ is a diffeomorphism\footnote{In this case, the Hamiltonian function $H$ is called \emph{hyperregular} \cite[page 73]{Libermann_marle87}}.
\item[(b)] $S(x)=e^\top x$ with some $e\in\mathbb{R}^n\backslash\{0\}$, i.e., the total entropy is linear in the state.
\end{enumerate}
\end{as}
\noindent This assumption may be interpreted in terms of thermodynamics in the following way. Item (a) 
is satisfied for simple thermodynamical systems by the strict convexity property of the internal energy function with respect to a set of independent extensive variables \cite{Callen85}. For complex systems, using the fact that the state variables are the extensive thermodynamic variables, the item (a) is equivalent to state that there is no constraint on the intensive variables, hence that the subsystems are not in thermodynamic equilibrium \footnote{If this was the case, one would have to consider more general formulations such as \cite{Entropy_2018,Eberard_RMP07}.}. Item (b) of the assumption corresponds the inclusion of entropy of each subsystem in the state vector and their balance equation in the dynamical equations.
\noindent Note that Lemma \ref{l:manifold} implies that under these assumptions the set of thermodynamic equilibria $\calT$ is a $C^1$-submanifold of $\R^n$.

We briefly verify Assumption~\ref{as:ham} for the two examples of the previous section.

\begin{ex}
{\bf (a)} We consider the heat exchanger from Example \ref{ex:heat}. Here, the state space is $\XX = \R^2$. The entropy $S$ is given by $S(x) = e^\top x$, where $e = \begin{bmatrix}1&1\end{bmatrix}^\top$, and the Hamiltonian $H$ is a $C^2$-function with
$$
H_x(S_1,S_2) = T_{\rm ref}\bvek{e^{(S_1-S_{\rm ref})/c_1}}{e^{(S_2-S_{\rm ref})/c_2}}.
$$
Hence, $\YY = (0,\infty)\times (0,\infty)$, and $H_x : \XX\to\YY$ is bijective with the continuously differentiable inverse function $H_x^{-1} : \YY\to\XX$,
$$
H_x^{-1}(T_1,T_2) = S_{\rm ref}\bvek 11 - (\ln T_{\rm ref})\bvek{c_1}{c_2} + \bvek{c_1\ln T_1}{c_2\ln T_2}.
$$
{\bf (b)} In case of the gas-piston system (Example \ref{ex:gp}), the entropy is a part of the state: $S(x) = \begin{bmatrix}1&0&0\end{bmatrix}x$. Concerning the energy variables $S$, $V$, and $p$, we allow $S$ and $p$ to assume any real value and $V$ only assumes values in an interval $(0,V_{\max})$, where $V_{\max}>0$ is the volume of the chamber. The state space is thus given by $\XX = \R\times (0,V_{\max})\times\R$. The image of $H_x$ then equals
$$
\YY = 
\left\{\bvek T{-P+\frac{mg}{A}} : T>0,\;P > \frac{RNT}{V_{\max}}\right\}\times\R.
$$
This easily follows from the representation
$$
C_0\cdot e^{\frac{2S}{3RN}}\bvek{(P_0V)^{-2/3}}{-RNP_0(P_0V)^{-5/3}} + \bvek 0{\frac{mg}A}
$$
of the first two components of $H_x(S,V,p)$, where
$$
C_0 = T_0^{5/3}(RN)^{2/3}e^{-\frac{2s_0}{3R}}.
$$
The inverse map $H_x^{-1} : \YY\to\XX$ is given by
$$
H_x^{-1}(T,-P+\tfrac{mg}A,v) = 
\begin{bmatrix}
Ns_0 + \frac{5RN}{2}\ln\frac{P_0}{P}\ln\frac T{T_0}\,,  &  \frac{RNT}{P}\,,  &  mv
\end{bmatrix}^\top,
$$
which is clearly continuously differentiable.
\end{ex} 

\begin{rem}
In other works treating the gas-piston system, see, e.g., \cite{Ramirez_EJC13,zla23}, 
an additional energy variable in the state, namely the position $z$ of the piston is considered. The corresponding co-energy variable is the constant gravitational force $F_g = mg$. Hence, in this extended set of coordinates, the gradient of the Hamiltonian $H_x$ is {\em not} a diffeomorphism. However, the weaker condition in Lemma \ref{l:manifold} is satisfied and implies that the set of thermodynamic equilibria is also a $C^1$-submanifold in this case.
\end{rem}

\begin{lem}
If $\linspan\{e\}\cap\YY\neq\emptyset$, the curve $H_x^{-1}(\linspan\{e\}\cap\YY)$ is contained in $\calT_{\rm opt}$ and thus, in particular, in $\calT$.
\end{lem}
\begin{proof}
If $\hat x\in H_x^{-1}(\linspan\{e\}\cap\YY)$, then we have $H_x(\hat x)\in\linspan\{e\}$ $= \linspan\{S_x(\hat x)\}$. Hence, $(\hat x,0)\in\calS$ by Proposition \ref{p:oss}, which implies $\hat x\in\calT_{\rm opt}$.
\end{proof}

\begin{ex}
In the case of the heat exchanger (Example \ref{ex:heat}),
$$
\linspan\{e\}\cap\YY = (Je)^\perp\cap\YY = \{(\la,\la)^\top : \la > 0\}.
$$
It is then easily seen that $\calT = \calT_{\rm opt} = H_x^{-1}((Je)^\perp)$ coicides with the affine subspace $v_0 + \linspan\{v_1\}$, where $v_0 = S_{\rm ref}\cdot e$ and $v_1 = [c_1,c_2]^\top$. In particular, $\calT = \calT_{\rm opt} = \linspan\{e\}$, if $c_1=c_2$.
\end{ex}

The proof of the following proposition can be found in the Appendix.


\begin{prop}\label{p:distance}
Let $V := \bigcap_{k=1}^N(J_ke)^\perp$ and assume that $V\cap\YY\neq\emptyset$. Then for each compact set $K\subset\XX$ we have
$$
\dist(x,\calT)^2\,\lesssim\,\sum_{k=1}^N\gamma_k(x,H_x)\{S,H\}_{J_k}^2(x),\quad x\in K.
$$
\end{prop}

Next, we introduce the manifold turnpike property that we will prove in the remainder of this section for state transition problems, and in the following section for output stabilization problems. This turnpike property resembles an integral version of the measure turnpike property introduced in \cite{tudo:faulwasser22e}.

\begin{defn}[Integral state manifold turnpike property]
Let $\ell\in C^1(\XX\times\mathbb U)$, $\vphi\in C^1(\XX)$, $\Phi\subset\XX$ a closed set and $f\in C^1(\R^{n+m})$. We say that a general OCP of the form
\begin{align}
\begin{split}\label{e:lin_OCP}
\min_{u\in L^1(0,t_f;\mathbb{U})}\,&\varphi(x(t_f)) + \int_0^{t_f} \ell(x(t),u(t))\,dt\\
\text{s.t. }\dot x = &f(x,u),\quad x(0)=x^0,\quad x(t_f)\in\Phi
\end{split}
\end{align}
has the {\em integral turnpike property} on a set $S_{\rm tp}\subset\XX$ with respect to a manifold $\mathcal{M}\subset\R^n$ if for all compact $X^0\subset S_{\rm tp}$ there is a constant $C>0$ such that for all $x^0\in X^0$ and all $t_f>0$, each optimal pair $(x^\star ,u^\star)$ of the OCP \eqref{e:lin_OCP} with initial datum $x^\star(0)=x^0$ satisfies
\begin{align}\label{e:integral_tp}
\int_0^{t_f}\dist^2\big(x^\star(t),\mathcal{M}\big)\,dt\,\le\,C.
\end{align}
\end{defn}

Due to the lack of uniform bounds on the Hessian $H_{xx}$ of the Hamiltonian,  Proposition~\ref{p:distance} holds only on compact sets. To apply this result to optimal trajectories and to render the involved constants uniformly in the horizon, we now assume that optimal trajectories are uniformly bounded in the horizon.

\begin{as}\label{as:comp}
For any compact set of initial values $X^0\subset\XX$ there is a compact set $K\subset\XX$ such that for all horizons $t_f$ each corresponding optimal state trajectory of \eqref{eq:phOCP} with initial datum $x^0\in X^0$ and horizon $t_f$ is contained in $K$. 
\end{as}

Define the following set of initial values which can be controlled to a state $\bar x\in\calT_{\rm opt}$ that can be further steered to the terminal set $\Phi$.
\begin{align*}
\mathcal{C}(\calT_{\rm opt},\Phi)
:= &\big\{x^0\in \mathbb{R}^n :\exists\,\bar{x}\in\calT_{\rm opt}\text{ s.t.}\\
&\hspace*{1.7cm}\exists\,t_1\geq 0,u_1 \in L^1(0,t_{1},\mathbb{U}) \text{ s.t. } x(t_1,u_1,x^0) = \bar x,\\
&\hspace*{1.7cm}\exists\,t_2\geq 0,u_2 \in L^1(0,t_{2},\mathbb{U}) \text{ s.t. } x(t_2,u_2,\bar{x})\in\Phi\big\}.
\end{align*}
We now provide a turnpike result for the optimal control problem~\eqref{eq:phOCP} subject to coupled reversible-irreversible dynamics~\eqref{eq:RIPHS}. A similar result was proven in the conference paper~\cite{MascPhil22} for the smaller class of irreversible systems, that is, $J_0\equiv 0$ and $N=1$ in \eqref{eq:RIPHS}.

\begin{thm}\label{thm:turnpike2}
Let Assumption~\rmref{as:comp} hold and furthermore assume that $\calC(\calT_{\rm opt},\Phi)\neq\emptyset$. Then, for any pair $\alpha_1,\alpha_2 > 0$ the OCP \eqref{eq:phOCP} has the integral turnpike property on the set $S_{\rm tp} = \calC(\calT_{\rm opt},\Phi)$ with respect to $\calT$.
\end{thm}
\begin{proof}[Proof \braces{Sketch}.]
The proof follows the lines of the proof of Theorem~8 in \cite{MascPhil22}. Let $X^0\subset S_{\text{tp}}$ be compact, let $x^0\in X^0$, and let $(x^*,u^*)$ be an optimal pair for \eqref{eq:phOCP}. By optimality, any control $u\in L^1(0,t_f;\mathbb{U})$, which is feasible for \eqref{eq:phOCP} with corresponding state trajectory $x= x(\,\cdot\,;u,x^0)$ satisfies
\begin{align*}
\int_0^{t_f} \ell_{\alpha_1,\alpha_2}(x^*(t),u^*(t))\,\text{d}t \leq  \int_0^{t_f} \ell_{\alpha_1,\alpha_2}(x(t),u(t))\,\text{d}t.
\end{align*}
Choose times $t_1,t_2$ and controls $u_1,u_2$ as in the definition of $\calC(\calT_{\rm opt},\Phi)$. For $t_f\le t_1+t_2$, it follows from Assumption \ref{as:comp} that \eqref{e:integral_tp} is satisfied with $C = (t_1+t_2)\sup\{\dist(x,\calM) : x\in K\}$. Let $t_f > t_1+t_2$ and define the constructed control 
\begin{align*}
u(t):=\begin{cases}
u_1(t) \qquad &t\in [0,t_1]\\
\bar{u} \qquad& t \in (t_1,t_f-t_2)\\
u_2(t-(t_f-t_2)) \qquad &t\in [t_f-t_2,t_f]
\end{cases}.
\end{align*}
This control steers the system from $x^0$ via $\bar x\in\calT_{\rm opt}$ (where it remains from time $t_1$ until $t_f - t_2$) to the terminal region $\Phi$. The middle part $\bar{u}$ is the steady state control that is required to stay at the controlled equilibrium $\bar{x}$ with zero stage cost, i.e., $(\bar x,\bar u)\in\calS$. Hence, we have
\begin{align*}
\int_0^{t_f} \ell_{\alpha_1,\alpha_2}(x(t),u(t))\,\text{d}t
&= \left(\int_0^{t_1} + \int_{t_f-t_2}^{t_f}\right) \ell_{\alpha_1,\alpha_2}(x(t),u(t))\,\text{d}t = \sum_{i=1}^2\int_0^{t_i}\ell_{\alpha_1,\alpha_2}(x_i(t),u_i(t))\,\text{d}t,
\end{align*}
where $x_i$ is the state response to $u_i$ on $[0,t_i]$. Note that this expression is in fact independent of $t_f$. Denote its norm by $C_1$. Making use of Proposition \ref{p:distance} and \eqref{eq:energybalance}, we obtain
\begin{align*}
\int_0^{t_f}\dist^2(x^\star,\calT)\,dt
&\le C_2\int_0^{t_f}\sum_{k=1}^N\gamma_k(x^\star,H_x(x^\star))\{S,H\}_{J_k}^2(x^\star)\,dt\\
&\le\tfrac{C_2}{\alpha_2T_0}\left(C_1 - \alpha_1[H(x^\star(t_f)) - H(x^0)]\right) + C_2 [S(x^\star(t_f)) - S(x^0)].
\end{align*}
Since $H$ and $S$ are continuous, the right-hand side can be estimated independently of $t_f$ due to Assumption \ref{as:comp}.
\end{proof}

\begin{rem}[Relaxing the conditions of Theorem~\ref{thm:turnpike2}]~\\
In the proof of Theorem~\ref{thm:turnpike2}, reachability of a controlled equilibrium is used to bound the cost functional uniformly with respect to the time horizon $t_f$. This reachability can be relaxed to exponential reachability of the subspace, i.e., for all $x^0$ there is a measurable control $u:[0,\infty)\to \mathbb{U}$ such that
\begin{align*}
    \dist(x(t;x^0,u),\bar{x})\leq Me^{-\omega t}\|x^0\|
\end{align*}
with $M \geq 1$ and $\omega > 0$, cf.\ e.g.\ \cite[Assumption 1]{FaulKord17}. Correspondingly, one has to impose reachability of the terminal state from a ball around the steady state with radius depending on the time horizon and the compact set of initial values. Other decay rates are also possible, as long as they are integrable on the positive real line. 
\end{rem}

%
%

\section{Optimal output stabilization}\label{sec:outputstab}
\noindent Let a horizon $t_f\geq 0$, an initial state $x^0\in \mathbb{R}^n$ be given and assume that $\mathbb{U}\subset\mathbb{R}^m$ is compact and convex. Consider an output matrix $C\in \R^{p\times n}$ and $y_\text{ref}\in\operatorname{im}(C)$. Let $\alpha_1,\alpha_2 \in \R^+$ and consider the optimal control problem with  
\begin{align}
\label{eq:phOCP_output}\tag{phOCP-stabilization}
\begin{split}
\min_{u\in L^1(0,t_f;\mathbb{U})} \int_0^{t_f}& \big[\|Cx(t) - y_\text{ref}\|^2 + \ell_{\alpha_1,\alpha_2}(x(t),u(t))\big]\,\text{d}t \\
\text{s.t. } \eqref{eq:RIPHS},\quad &x(0)=x^0
\end{split}
\end{align}
Compared to \eqref{eq:phOCP}, we do not consider a state transition with a terminal set but aim to stabilize an output $y_\text{ref}$ in the cost functional. Again, we use the energy balance to reformulate the OCP, where $\wt\ell_{\alpha_1,\alpha_2}$ is the stage cost in \eqref{eq:phOCP_output}:
\begin{align}
\begin{split}
\int_0^{t_f} \wt\ell_{\alpha_1,\alpha_2}(x,u)\,\text{d}t= & \alpha_1\left[H(x(t_f))-H(x(0))\right] + \int_0^{t_f}\|Cx(t) - y_\text{ref}\|^2\,\text{d}t\\ &\quad+ \alpha_2T_0\left(S(x(0))- S(x(t_f)) + \int_0^{t_f}\gamma(x,H_x)\{S,H\}^2_J\,\text{d}t\right).
\end{split}
\end{align}
The long-term behavior is now governed by two terms. As in \eqref{eq:energybalance}, the term $\int_0^{t_f} \gamma(x,H_x) \{S,H\}^2_J\,\text{d}t$ corresponds to the distance to the set of thermodynamic equilibria, as shown in Subsection~\ref{subsec:statetrans_turnpike}. However, we obtained the additional output stabilization term $\int_0^{t_f}\|Cx(t) - y_\text{ref}\|^2\,\text{d}t$ penalizing the distance to the preimage of $y$ under $C$ as shown in the following.
\begin{lem}\label{l:distance_ker}
Let $y\in\operatorname{im}(C)$. Then
\begin{align*}
\|Cx-y\|\,\asymp\,\dist(x,C^{-1}\{y\}),\quad x\in\R^n,
\end{align*}
where $C^{-1}\{y\}$ is the preimage of $y$ under $C$.
\end{lem}
\begin{proof}
Let $A = C^\top C$, and let $x_0\in\R^n$ be such that $Cx_0 = y$. By \cite[Lemma 3.15]{Faulwasser2021} we have $z^\top Az\asymp\dist^2(z,\ker A)$ for $z\in\R^n$. Setting $z =  x - x_0$, we note that
$$
\dist(z,\ker A) = \dist(x,x_0+\ker C) = \dist(x,C^{-1}\{y\})
$$
and $z^\top Az = \|C(x-x_0)\|^2 = \|Cx - y\|^2$.
\end{proof}
\noindent We now provide the turnpike result for the output stabilization problem \eqref{eq:phOCP_output}.

\begin{thm}[Turnpike for output  stabilization]\label{thm:turnpike_tracking}
Let Assumption \rmref{as:comp} hold for the optimal control problem \eqref{eq:phOCP_output}, and suppose that the set $\calC(\calT_{\rm opt}',\R^n)$ is non-empty, where $\calT_{\rm opt}' := \calT_{\rm opt}\cap C^{-1}\{y_{\rm ref}\}$. Then, for any pair $\alpha_1,\alpha_2 > 0$ the OCP \eqref{eq:phOCP_output} has the integral turnpike property on the set $S_\text{tp} = \calC(\calT_{\rm opt}',\R^n)$ with respect to $\calT$ and $ C^{-1}\{y_{\rm ref}\}$. If $\calT$ is a subspace, then this turnpike property holds with respect to $\calT\cap C^{-1}\{y_{\rm ref}\}$.
\end{thm}
\begin{proof}
The proof follows analogously to the proof of Theorem~\ref{thm:turnpike2}. As we do not have to fulfill a terminal constraint, we can consider the constructed control
\begin{align*}
u(t):=\begin{cases}
u_1(t) \qquad &t\in [0,t_1]\\
\bar{u} \qquad& t \in (t_1,t_f].
\end{cases}
\end{align*}
where $u_1$ steers $x^0\in \calC(\calT_{\rm opt}',\R^n)$ into $\calT_{\rm opt}'$. This construction allows (analogously to the argumentation of the proof of Theorem~\ref{thm:turnpike2}) to bound the cost of the constructed trajectory by a constant $C_1\geq 0$ independent of the horizon $t_f>0$. Then, an optimality argument and invoking Assumption~\ref{as:comp} yields
\begin{align*}
    \int_{t_0}^{t_f} \|Cx^*(t) - y_\mathrm{ref}\|^2 + \sum_{k=1}^N\gamma_k(x^\star,H_x(x^\star))\{S,H\}_{J_k}^2(x^\star)\,\mathrm{d}t \leq C_2
\end{align*}
with $C_2\geq 0$ independent of $t_f$. Proposition~\ref{p:distance} and Lemma~\ref{l:distance_ker} yield
\begin{align}\label{eq:ineqproof}
        \int_{t_0}^{t_f} \dist(x^*(t),C^{-1}\{y_\mathrm{ref}\}) + \dist(x^*(t),\calT)\,\mathrm{d}t \leq C_3
\end{align}
with $C_3\geq 0$ independent of $t_f$. This implies the claimed manifold turnpike property w.r.t.\ $C^{-1}\{y_\mathrm{ref}\}$ and $\calT$. 

If $\calT$ is a subspace, we can invoke Lemma~\ref{app:distances} such that for all $x\in \R^n$,
\begin{align*}
    \dist(x,C^{-1}\{y_\mathrm{ref}\} \cap \calT) \lesssim \dist(x^*(t),C^{-1}\{y_\mathrm{ref}\}) + \dist(x^*(t),\calT).
\end{align*}
Inserting this inequality into \eqref{eq:ineqproof} implies the manifold turnpike property w.r.t.\ the intersection $\calT\cap C^{-1}\{y_{\rm ref}\}$.
\end{proof}

\begin{ex}
We briefly comment on an application of this result to a heat exchanger network, which we will later also inspect numerically in Subsection~\ref{subsec:heat_network}. To this end, consider a heat exchanger network consisting of three identical subsystems. A straightforward extension of the modeling performed in Example~\ref{ex:heat} reveals that for a system of three coupled heat exchangers the thermodynamic equilibria are given by the two-dimensional subspaces
\begin{align*}
\mathcal{T}_1 = \{(S_1,S_2) : T_1=T_2\}   \qquad \text{and} \qquad
\mathcal{T}_2 = \{(S_2,S_3) : T_2=T_3\}.
\end{align*}
Note that, for simplicity we assumed that the constants in the temperature-entropy relation are identical for all subsystems. Clearly, we have the one-dimensional subspace $\mathcal{T} := \mathcal{T}_1 \cap\mathcal{T}_2 = \{S_1=S_2=S_3\}$. Let us briefly discuss three possible output configurations in view of the OCP \eqref{eq:phOCP_output}:

\begin{itemize}
\item[1)] No output in the cost functional: As the output term is not present, setting $\Phi = \mathbb{X}$ in Theorem~\ref{thm:turnpike2} yields a subspace turnpike towards the one-dimensional subspace~$\mathcal{T}$.
\item[2)] Output stabilization with scalar output and prescription of the temperature (or equivalently of the entropy) in subsystem $i\in \{1,2,3\}$. In this case, $C^{-1}\{y_\mathrm{ref}\} = \{S_i= y_\mathrm{ref}\}$, and we obtain the zero-dimensional turnpike manifold $\mathcal{T}\cap C^{-1}\{y_\mathrm{ref}\} = \{S_1=S_2=S_3=y_\mathrm{ref}\}$. Such an example will be considered in Subsection~\ref{subsec:heat_num2}.
\item[3)] Stabilizing temperature/en\-tropy in two subsystems: Here $C\in \R^{2\times 3}$ and $C^{-1}\{y_\mathrm{ref}\} = \{S_i = (y_\mathrm{ref})_i, S_j =(y_\mathrm{ref})_j \}$ for some $i,j\in \{1,2,3\}$, $i\neq j$ such that the intersection with $\mathcal{T}=\{S_1=S_2=S_3\}$ is empty if there is $i,j\in \{1,2,3\}$ such that $(y_\mathrm{ref})_i\neq (y_\mathrm{ref})_j$. In this case, Theorem~\ref{thm:turnpike_tracking} is not applicable. However, we will present a numerical example with a network of heat exchangers in Subsection~\ref{subsec:heat_network} and observe a turnpike-like behavior towards an optimal tradeoff.
\end{itemize}
\end{ex}

\begin{rem}
    We briefly comment on a combination of the results of Section~\ref{sec:statetrans} and Section~\ref{sec:outputstab} for problems involving tracking terms and terminal conditions. In this case, as the long-term behavior of optimal states is governed by the integral terms in the cost functional, a turnpike result similar to the case output tracking without terminal conditions is to be expected. More precisely, one could straightforwardly adapt Theorem~\ref{thm:turnpike_tracking} to terminal conditions $x(t_f)\in \Phi$ by replacing $\mathcal{C}(\mathcal{T}'_\mathrm{opt},\R^n)$ with $\mathcal{C}(\mathcal{T}'_\mathrm{opt},\Phi)$ and including also the leaving arc in the constructed control as in the proof of Theorem~\ref{thm:turnpike2}.
\end{rem}

\section{Numerical results}\label{sec:numerics}
\noindent In this part, we conclude an extensive numerical case study to illustrate the results of Sections~\ref{sec:statetrans} and~\ref{sec:outputstab}, in particular the turnpike property proven in Theorems~\ref{thm:turnpike2} for state-transition problems and in Theorem~\ref{thm:turnpike_tracking} for output stabilization. For time discretization, we use an implicit midpoint rule with time step $\Delta t = 0.01$, discretize the cost functional with the rectangular rule and solve the corresponding optimal control problem with CasADi~\cite{AndeGill2019}.

\subsection{Stabilization of a heat exchanger}\label{subsec:heat_num2}
\begin{figure}[ht]
\centering
\includegraphics[width=\linewidth]{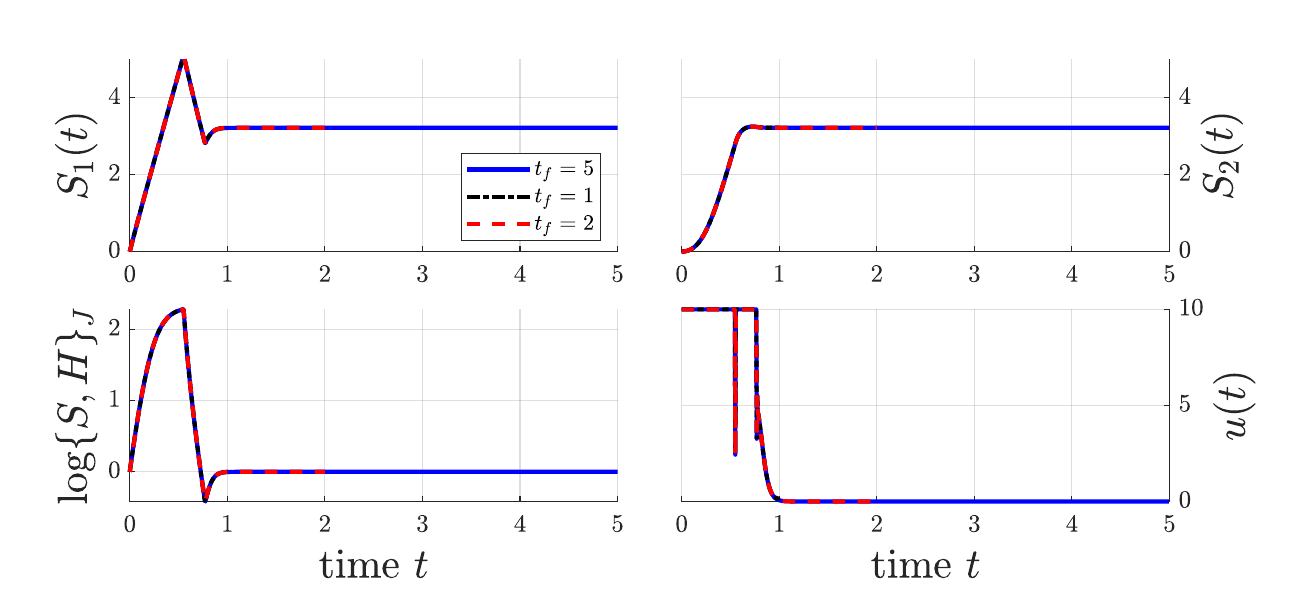}
\caption{Depiction of the optimal intensive variable, entropy bracket and control over time for output stabilization.}
\label{fig:heat_tracking}
\end{figure}

\noindent 
First, we provide a numerical example with a heat exchanger, cf.\ Example~\ref{ex:heat}. While a set point transition task was considered in \cite{MascPhil22}, we now aim to track a desired reference temperature of $T_\mathrm{ref}=25$ degrees.
For the parameters of the cost functional in \eqref{eq:phOCP_output}, we choose $\alpha_1=0$ and $T_0\alpha_2=1$, i.e., we minimize the entropy production. To obtain a linear output of the system state, we translate the aim temperature into an entropy, that is, we aim to minimize, in addition to the entropy creation, the term
\begin{align*}
|S_2(t) - S_\mathrm{ref}|^2 =  |Cx(t) - S_\mathrm{ref}|^2
\end{align*}
with $S_\mathrm{ref}=\log(25)$ and $C=\big[0\quad 1\big]$. Correspondingly, we have $C^{-1}\{S_{\rm ref}\} = \{(S_1,S_2)\in \mathbb{R}^2 : S_2=S_\mathrm{ref}\}$ which is a one-dimen\-sional affine subspace of $\R^2$. In addition, the manifold of thermodynamic equilibria is $\mathcal{T}=\{(S_1,S_2)\in \R^2:S_1=S_2\}$. Thus, Theorem~\ref{thm:turnpike_tracking} yields a turnpike w.r.t.\
$$ C^{-1}\{y_{\rm ref}\}\cap \mathcal{T} = \{(S_1,S_2)\in \R^2 : S_1=S_2=S_\mathrm{ref}\}.
$$
This can be observed in Figure~\ref{fig:heat_tracking}. There, both entropies approach the target entropy $S_\mathrm{ref}=\log(25)\approx 3.22$. In order to approach this value, the entropy (or equivalently, the temperature) in the first compartment has to be increased. This yields an increase of the entropy (or temperature) in the second compartment due to the heat flux across the wall, inevitably coupled to entropy generation, cf.\ bottom left of Figure~\ref{fig:heat_tracking}. Having reached the desired target entropy in the second compartment, the control is switched off such that the system is at equilibrium with $u=0$ and no entropy is generated. Thus, both integrands in the cost functional (approximately) vanish for this state: the output term in the cost vanishes as $S_2 \approx S_\mathrm{ref}$ and the entropy production vanishes as $S_1=S_2$ and thus $\{S,H\}_{J} = T_1 - T_2 \approx 0$.

\subsection{Set-point transition for a gas-piston system}
\noindent  Second, we consider the gas-piston system of Example~\ref{ex:gp}.
As an initial configuration of intensive and extensive variables in the thermodynamic domain, we choose
\begin{align*}
 S(0) = Ns_0,\quad V(0) =  \frac{NRT_0}{P_0}
\end{align*}
such that $\beta(S(0),P(0))=0$ as defined in \eqref{eq:InternalEnergyIdealGasAdd} and $T(0)=T_0$,  $P(0)=P_0$.
Moreover, we assume that the mechanical subsystem is at equilibrium, that is $v(0)=0$ and $p(0)=0$. 
The mass is also chosen such that the momentum ODE (last line of \eqref{dynamic_gas-piston}) is in equilibrium for zero velocity, i.e., 
$$
m = \frac{AP_0}{g} = 5.1644\;\text{kg}.
$$
This yields a controlled equilibrium of the system \eqref{dynamic_gas-piston} endowed with the input map \eqref{eq:NonlinearControlTransformation} when choosing the external temperature $u = T_0$. We summarize the chosen parameters in Table~\ref{tab:gaspiston}.

\begin{figure}[htb]
\centering
\includegraphics[width=1.\linewidth]{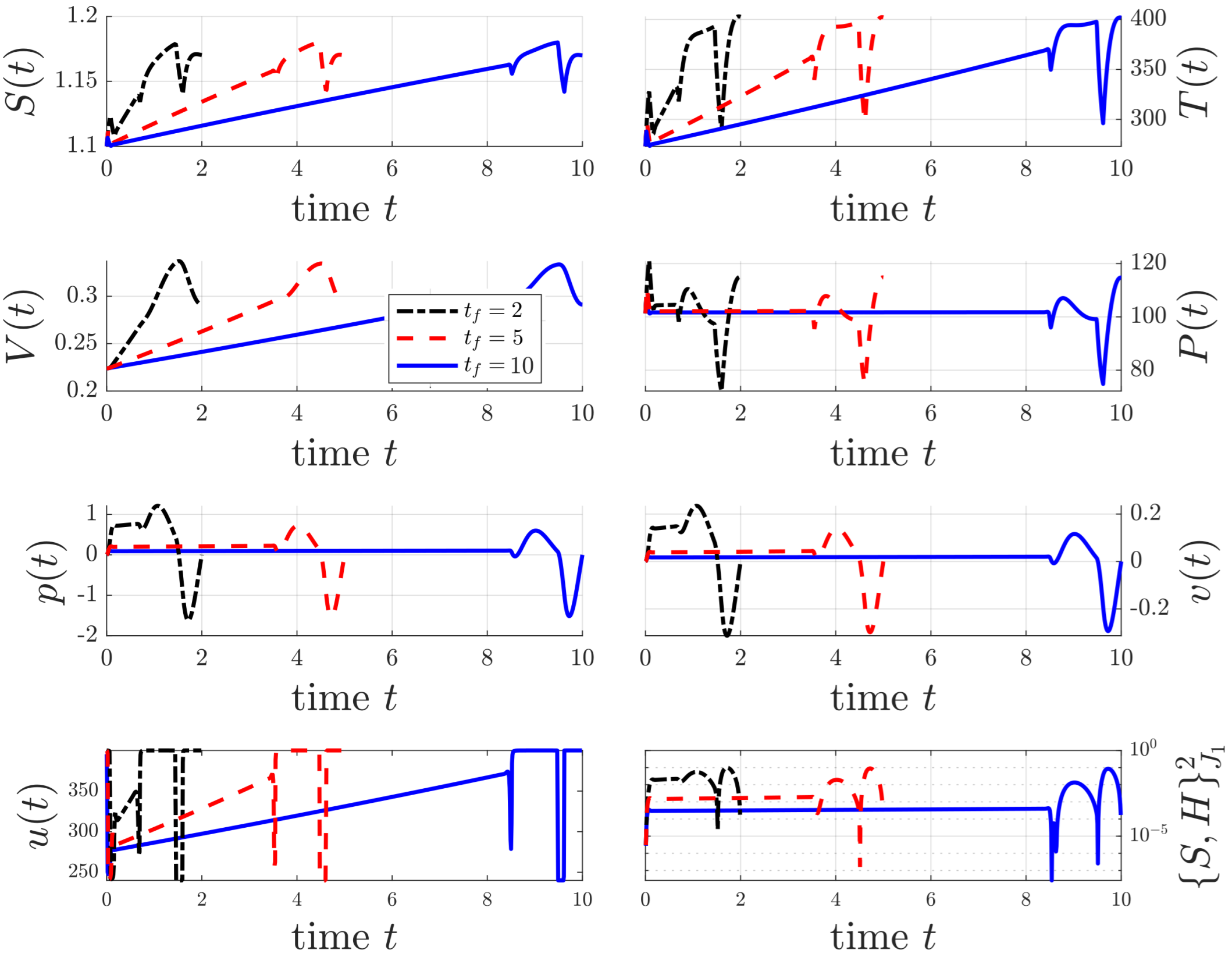}
\caption{Depiction of the optimal variables for the gas-piston system.}
\label{fig:gaspiston}
\end{figure}

We now aim to steer the piston to the target volume $V(t_f) = 1.3V(0)$ with target zero momentum $p(t_f) = 0$. In view of \eqref{eq:energybalance}, the choice of $\alpha_1$ does not matter due to fixed initial and terminal state, such that we set $\alpha_1=0$ and $T_0\alpha_2=1$. The optimal intensive variables, the corresponding extensive quantities, the entropy production by means of the Poisson bracket and the optimal control are depicted in Figure~\ref{fig:gaspiston}.

\noindent To change the volume of the piston, we need to induce a temperature increase by means of the heating jacket. As a result, the volume increases, whereas the pressure stays constant.
In view of Theorem~\ref{thm:turnpike2}, we observe a turnpike in the velocity component, as $\{S,E\}_{J_1} = v$, that is, the manifold of thermodynamic equilibria consists of the states with zero velocity. The longer the time horizon, the smaller the velocity for the majority of the time.

\begin{center}
\begin{table}[ht]\label{t:ps}
\begin{tabular*}{\columnwidth}{ll|l}
\hline 
\hline
Number of moles $N$ &[$-$]&0.01\\
Reference pressure $P_0$ &[Pa]&101.325 \\
Reference Temperature $T_0$ &[K]&273\\
Reference molar entropy $s_0$ &[kJ K$^{-1}$ mol$^{-1}$]&0.1100\\
Mass $m$ &[kg]& 5.1644\\
Friction coefficient $\nu$ &[$-$]&10\\
Area $A$&[m$^2$]&0.5\\
heat conduction coeff.\ $\lambda_e$&[$-$]&1\\	
\hline 
\hline
\end{tabular*}
\caption{Parameters and reference values}
\label{tab:gaspiston}
\end{table}
\end{center}

\subsection{Network of heat exchangers}\label{subsec:heat_network}
\begin{figure}[htb]
	\resizebox{.8\linewidth}{!}{%
	\begin{tikzpicture}[
	node distance= 9mm and 0mm,
	rect/.style = {draw=black, fill=gray!20, 
		minimum width=15pt, minimum height = 10pt},
	every edge/.style = {draw, -Latex},
	node distance=2.5cm
	]
	\node[rect] (a) {$S_1,T_1$};
	\node[rect, right of= a] (b) {$S_2,T_2$};
    \node[rect, right of = b] (c) {$S_3,T_3$};
	\node[rect,  right of = c] (d) {$S_4,T_4$};
	\node[rect] (e) [below=.6cm of c]{$S_5,T_5$};	
 	\node[] (u1) [above=.6cm of a] {};
 	\node[] (u2) [above=.6cm of d] {};
 	\node[] (u3) [left=.6cm of e] {};

	\path (u1) edge["{$u_1$}"] (a.north)
    (u2) edge["{$u_2$}"] (d.north)
    (u3) edge["{$u_3$}"] (e.west)
    (a) edge["\scriptsize{$\lambda_1(T_1-T_2)$}" swap] (b)
	(b) edge["\scriptsize{$\lambda_2(T_2-T_3)$}"] (c)
	(c) edge["\scriptsize{$\lambda_3(T_3-T_4)$}"] (d)
 	(c.south) edge["\scriptsize{$\lambda_4(T_3-T_5)$}"] (e.north);
 	\end{tikzpicture}}
	\caption{Network of five compartments exchanging heat.}
\label{fig:coupledheat}
\end{figure}
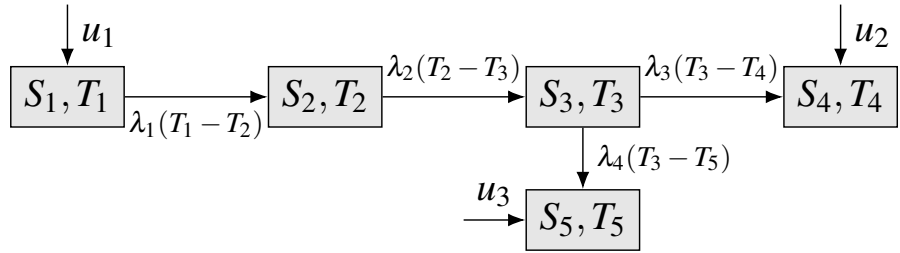

\begin{figure}[htb]
    \centering\includegraphics[width=\linewidth]{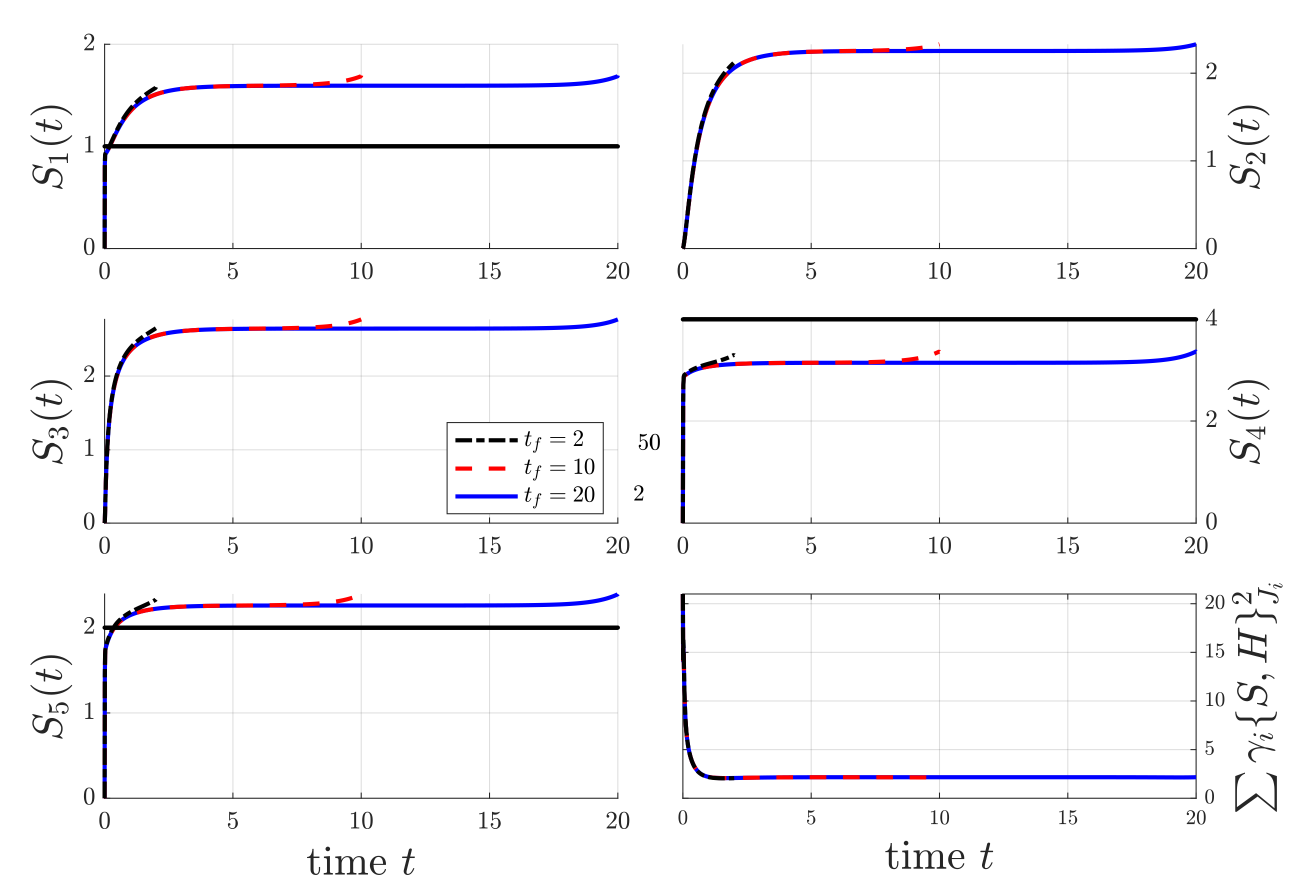}
    \caption{Entropies in the individual compartments for the network of heat exchangers and generated entropy.}
    \label{fig:coupledheat_results}
\end{figure}

\noindent Last, we present an example with five compartments exchanging heat, which are coupled as depicted in Figure~\ref{fig:coupledheat}.
The corresponding state variables are the entropies in the compartments $S_i$, $i=1,\ldots,5$ with Hamiltonian energies $H_i(S_i) = e^{S_i}$. $i=1,\ldots,5$. Along the lines of the Example~\ref{ex:heat} 
we get by Fourier's law and continuity of the heat flux for the first, fourth and fifth compartment
\begin{align*}
\lambda_1(T_1-T_2) &= -\frac{\mathrm{d}}{\mathrm{d}t} H_1(S_1(t)) = -T_1\frac{\mathrm{d}}{\mathrm{d}t}{S}_1(t),\\
\lambda_2(T_3-T_4) &= \phantom{-}\frac{\mathrm{d}}{\mathrm{d}t} H_4(S_4(t)) = \phantom{-}T_4\frac{\mathrm{d}}{\mathrm{d}t}{S}_4(t),\\
\lambda_3(T_3-T_5) &= \phantom{-}\frac{\mathrm{d}}{\mathrm{d}t} H_5(S_5(t)) = \phantom{-}T_5\frac{\mathrm{d}}{\mathrm{d}t}{S}_5(t).
\end{align*}
For the second compartment, we compute
\begin{align*}
T_2\frac{\mathrm{d}}{\mathrm{d}t}{S}_2(t) &= \frac{\mathrm{d}}{\mathrm{d}t} H_2(S_2(t)) \\&= \lambda_1(T_1-T_2)-\lambda_2(T_2-T_3),
\end{align*}
and for the third compartment,
\begin{align*}
T_3\frac{\mathrm{d}}{\mathrm{d}t}{S}_3(t) &= \frac{\mathrm{d}}{\mathrm{d}t} H_3(S_3(t)) \\&= \lambda_2(T_2-T_3)-\lambda_2(T_3-T_4)- \lambda_4(T_3-T_5).
\end{align*}
Thus, corresponding to the four coupling interfaces between the five compartments we define the skew-symmetric structure matrices
\begin{align*}
J_1 &= \begin{bmatrix}
0&-1&0&0&0\\
1&\phantom{-}0&0&0&0\\
0&\phantom{-}0&0&0&0\\
0&\phantom{-}0&0&0&0\\
0&\phantom{-}0&0&0&0
\end{bmatrix},\quad 
J_2= \begin{bmatrix}
0&0&\phantom{-}0&0&0\\
0&0&-1&0&0\\
0&1&\phantom{-}0&0&0\\
0&0&\phantom{-}0&0&0\\
0&0&\phantom{-}0&0&0
\end{bmatrix},\\
J_3 &=  \begin{bmatrix}
0&0&0&\phantom{-}0&0\\
0&0&0&\phantom{-}0&0\\
0&0&0&-1&0\\
0&0&1&\phantom{-}0&0\\
0&0&0&\phantom{-}0&0
\end{bmatrix},\quad 
J_4 =  \begin{bmatrix}
0&0&0&0&\phantom{-}0\\
0&0&0&0&\phantom{-}0\\
0&0&0&0&-1\\
0&0&0&0&\phantom{-}0\\
0&0&1&0&\phantom{-}0
\end{bmatrix}
\end{align*}
and the positive functions
\begin{align*}
\gamma_1(x,H_x) = \frac{\lambda_1}{T_1T_2},\qquad \gamma_2(x,H_x) = \frac{\lambda_2}{T_2T_3},\\
\gamma_3(x,H_x) = \frac{\lambda_3}{T_3T_4},\qquad \gamma_4(x,H_x) = \frac{\lambda_3}{T_3T_5}.
\end{align*}
The corresponding Poisson brackets giving rise to the irrever\-sible phenomena read
\begin{align*}
\{S,H\}_{J_1} = T_1-T_2,\quad \{S,H\}_{J_2} = T_2-T_3,\\
\{S,H\}_{J_3} = T_3-T_4, \quad \{S,H\}_{J_4} = T_3-T_5.
\end{align*}
such that
\begin{align*}
\frac{\mathrm{d}}{\mathrm{d}t}S_1(t) &= 
-\gamma_1(x,H_x)\{S,H\}_{J_1} T_2\\
\frac{\mathrm{d}}{\mathrm{d}t}S_4(t) &= \phantom{-}\gamma_3(x,H_x)\{S,H\}_{J_3} T_3\\
\frac{\mathrm{d}}{\mathrm{d}t}S_5(t) &= \phantom{-}\gamma_4(x,H_x)\{S,H\}_{J_4} T_3
\end{align*}
and
\begin{align*}
\frac{\mathrm{d}}{\mathrm{d}t}S_2(t)& = \gamma_1(x,H_x)\{S,H\}_{J_2}T_1-\gamma_2(x,H_x)\{S,H\}_{J_2} T_3\\
\frac{\mathrm{d}}{\mathrm{d}t}S_3(t)& = \gamma_2(x,H_x)\{S,H\}_{J_2}T_2-\gamma_3(x,H_x)\{S,H\}_{J_3} T_4\\&\qquad \qquad \qquad \qquad \quad-\gamma_4(x,H_x)\{S,H\}_{J_4} T_5.
\end{align*}
Eventually, we obtain with $x=(S_1,S_2,S_3,S_4,S_5)$ the dynamics
\begin{align*}
\frac{\mathrm{d}}{\mathrm{d}t}x(t) = \bigg(\sum_{i=1}^{4} \gamma_i(x,H_x)\{S,H\}_{J_i}J_i\bigg)H_x(x(t)).
\end{align*}
We endow our system with an entropy flow control at the first, the fourth and the fifth compartment and aim to track the temperature in these compartments. More precisely, we add the input term
\begin{align*}
g(x,H_x)u=
\begin{bmatrix}
1&0&0\\
0&0&0\\
0&0&0\\
0&1&0\\
0&0&1
\end{bmatrix}
\begin{bmatrix}
u_1\\u_2\\u_3
\end{bmatrix}
\end{align*}
and define the output stabilization term
\begin{align*}
\|Cx-y_\mathrm{ref}\|^2 =
\left\| \begin{bmatrix}
1&0&0&0&0\\
0&0&0&1&0\\
0&0&0&0&1
\end{bmatrix}x - \begin{bmatrix}
S_{1,\mathrm{ref}}\\
S_{4,\mathrm{ref}}\\
S_{5,\mathrm{ref}}
\end{bmatrix} \right\|^2.
\end{align*}
for the cost functional. If $S_{1,\mathrm{ref}}\neq S_{4,\mathrm{ref}}$ or 
$S_{1,\mathrm{ref}}\neq
S_{5,\mathrm{ref}}$, then Theorem~\ref{thm:turnpike_tracking} does not apply as
\begin{align*}
    &C^{-1}\left\{\begin{pmatrix}
S_{1,\mathrm{ref}},
S_{4,\mathrm{ref}},
S_{5,\mathrm{ref}}
\end{pmatrix}^\top\right\} = \left\{x\in \R^5 : x_1 = S_{1,\mathrm{ref}},\,x_4 = S_{4,\mathrm{ref}},\, x_5=S_{5,\mathrm{ref}}\right\}
\end{align*}
and
\begin{align*}
    \mathcal{T} = \left\{x\in \R^5 : x_1=x_2=\ldots=x_5\right\}
\end{align*}
such that
\begin{align}\label{eq:emptyintersec}
\begin{split}
       &C^{-1}\left\{\begin{pmatrix}
S_{1,\mathrm{ref}},
S_{4,\mathrm{ref}},
S_{5,\mathrm{ref}}
\end{pmatrix}^\top\right\} \cap \mathcal{T}_\mathrm{opt} \subset
       C^{-1}\left\{\begin{pmatrix}
S_{1,\mathrm{ref}},
S_{4,\mathrm{ref}},
S_{5,\mathrm{ref}}
\end{pmatrix}^\top\right\} \cap \mathcal{T} = \emptyset.
\end{split}
\end{align}
The results are shown in Figure~\ref{fig:coupledheat_results}. We see that $S_1$, $S_4$ and $S_5$ do not approach the given reference value, as this value would not allow for zero entropy creation. We rather observe a trade-off between the terms in the cost functional: The temperature (or equivalently the entropy) in the fifth and the first compartment is higher than the given reference, whereas in the fourth compartment the temperature is lower as desired. As can be seen in the bottom right plot of Figure~\ref{fig:coupledheat_results}, also the entropy creation is not small for the majority of the time as we have a trade-off between output stabilization and entropy creation. The reason for this trade-off is that due to \eqref{eq:emptyintersec} there is no state for which both terms in the cost functional vanish, which in particular prohibits an application of Theorem~\ref{thm:turnpike_tracking}. However, we still can observe a turnpike behavior, which should be investigated in future work.

\section{Conclusions and future work}\label{sec:conclusions}

\noindent 
We have considered an optimal control problem, intrinsically defined as minimizing the energy supply to the system, 
the irreversible entropy creation, or a linear combination of both, the minimal exergy destruction. To this end, we have formulated the physical model of the system as a reversible-irreversible port-Hamiltonian system, which is defined via a quasi-Poisson bracket and two functions: the total energy acting as a generating function and the total entropy function. We have characterized optimal state-control pairs of the steady state problem in terms of the manifold of thermodynamic equilibria. For dynamic state-transition and output stabilization problems we have derived conditions, under which the optimal solutions of the dynamic problem reside close to the manifold of thermodynamic equilibria for the majority of the time. Last, we have illustrated our results by means of various examples, including purely irreversible systems such as a network of heat exchangers or a reversible-irreversible gas-piston problem consisting of coupled mechanical and thermodynamic systems.

A topic for future work is the derivation of sufficient conditions to easily verify Assumption~\ref{as:comp}, e.g., based on a tailored detectability notion. Then, one may leverage the interplay between optimality and detectability to show the imposed uniform boundedness for all optimal trajectories emanating from a given compact set.
Moreover, the results of this paper will be extended in various ways. First, the results may be extended to a much wider class of systems including metriplectic and irreversible Hamiltonian systems as considered in \cite{goreac2024generating} \cite{Morrison_PhysRevE_2024}. 
Secondly, one may consider a different cost function including additionally the irreversible entropy production at the interfaces where a dissipative phenomenon takes place, as it has been defined in \cite{Maschke_submIFAC_WC23b}. 
Thirdly, these results might be adapted to the infinite-dimensional irreversible Port Hamiltonian systems as suggested in \cite{Ramirez_ChemEngSci_2022}, having in view applications such as the plug flow reactor and extending the stationary optimal control problem as in \cite{Johannessen_Energy_04_OptContrEntropy} to the dynamic set-point optimal control as presented here.

\bibliographystyle{abbrv}
\bibliography{irrevPH.bib}

\appendix
\section{Distances}\label{app:distances}
For a closed convex set $M\subset\R^n$ we denote by $P_M$ the orthogonal projection onto $M$.

\begin{lem}\label{l:subspaces}
Let $V_1,\ldots,V_N\subset\R^n$ be linear subspaces. Then 
\begin{align}\label{e:dists}
\dist\left(x,\bigcap_{k=1}^N V_k\right) \,\asymp\, \sum_{k=1}^N \dist(x,V_k)\quad\text{for $x\in\R^n$}.
\end{align}
\end{lem}
\begin{proof}
Since $\dist(x,V_j)\le\dist(x,\bigcap_{k=1}^NV_k)$ for $j=1,\ldots,N$, it is evident that the left-hand side of \eqref{e:dists} is not smaller than $\frac 1N\sum_{k=1}^N\dist(x,V_k)$. Moreover, it obviously suffices to prove the opposite inequality for $N=2$. For this, we set $V = V_1$ and $W = V_2$. Note that for any subspace $U$ we have $\dist(x,U) = \|P_{U^\perp}x\|$.

Suppose that the claim is false. Then there exists a sequence of vectors $(x_\ell)$ such that for all $\ell\in\N$
$$
\|P_{(V\cap W)^\perp}x_\ell\| \,>\, \ell\big(\|P_{V^\perp}x_\ell\| + \|P_{W^\perp}x_\ell\|\big).
$$
This in particular implies that $y_\ell := P_{(V\cap W)^\perp}x_\ell\neq 0$. Note that $P_{V^\perp}y_\ell = P_{V^\perp}x_\ell$ as $V^\perp\subset (V\cap W)^\perp$. Setting $z_\ell := y_\ell/\|y_\ell\|$, it follows that $\|P_{V^\perp}z_\ell\| + \|P_{W^\perp}z_\ell\| < \tfrac 1\ell$ for $\ell\in\N$. As $\|z_\ell\|=1$ for all $\ell\in\N$, we may assume that $z_\ell\to z$ as $\ell\to\infty$ for some vector $z$ with $\|z\|=1$. The latter inequality then shows that $P_{V^\perp}z = P_{W^\perp}z = 0$ and thus $z\in V\cap W$. But $z_\ell\in (V\cap W)^\perp$ for $\ell\in\N$, which implies $z=0$, contradicting $\|z\|=1$.
\end{proof}

\begin{lem}\label{l:app2}
Let $V\subset\R^n$ be a linear subspace such that $V\cap\YY\ne\emptyset$. Then we have
$$
\dist(x,H_x^{-1}(V\cap\YY))\,\lesssim\,\dist(H_x(x),V),\quad x\in K,
$$
for any compact set $K\subset\XX$.
\end{lem}
\begin{proof}
Suppose that the claim is false. Then there exist a compact set $K\subset\XX$ and a sequence $(x_n)\subset K$ such that
$$
\dist(x_n,H_x^{-1}(V\cap\YY)) > n\cdot\dist(H_x(x_n),V),\quad n\in\N.
$$
Since $K$ is compact, we may assume WLOG that $x_n\to x$ as $n\to\infty$ with some $x\in K$. Set $y_n := H_x(x_n)$ and $y := H_x(x)$. Choose $\veps>0$ such that $\ol{B_\veps(y)}\subset\YY$. We have
$$
\|P_{V^\perp}y_n\| = \dist(y_n,V) < \tfrac 1n\cdot \dist(x_n,H_x^{-1}(V\cap\YY))\to 0
$$
and thus also
$$
\|y-P_Vy_n\|\le\|y-y_n\| + \|P_{V^\perp}y_n\|\,\to\,0
$$
as $n\to\infty$. Therefore, there exists $n_0\in\N$ such that $y_n,P_Vy_n\in B_\veps(y)$ for $n\ge n_0$. In particular, for $n\ge n_0$ we have $P_Vy_n\in V\cap\YY$ and hence
\begin{align*}
n\|P_{V^\perp}y_n\|
&< \dist(H_x^{-1}(y_n),H_x^{-1}(V\cap\YY))
\le \|H_x^{-1}(y_n) - H_x^{-1}(P_Vy_n)\|\\
&\le\Big[\sup_{\xi\in [y_n,P_Vy_n]}\|D_\xi H_x^{-1}\|\Big]\cdot\|y_n - P_Vy_n\|
\le\Big[\sup_{\xi\in\ol{B_\veps(y)}}\|D_\xi H_x^{-1}\|\Big]\cdot\|P_{V^\perp}y_n\|,
\end{align*}
which is a contradiction for large $n$.
\end{proof}

\begin{proof}[Proof of Proposition \ref{p:distance}]
First of all, note that $\{S,H\}_{J_k}(x) = -H_x(x)^\top J_ke$ for $k=1,\ldots,N$. In particular, $\calT = H_x^{-1}(V\cap\YY)$. Also, $|w^\top v| = \|v\|\dist(w,v^\perp)$ for $v,w\in\R^n$. Hence, by Lemmas \ref{l:subspaces} and \ref{l:app2},
\begin{align*}
\dist(x,\calT)
&= \dist(x,H_x^{-1}(V\cap\YY))\,\lesssim\,\dist(H_x(x),V)\\
&\lesssim \sum_{k=1}^N\dist(H_x(x),(J_ke)^\perp)
\lesssim \sum_{k=1}^N\|J_ke\|\cdot\dist(H_x(x),(J_ke)^\perp)\\
&= \sum_{k=1}^N|H_x(x)^\top J_ke| = \sum_{k=1}^N|\{S,H\}_{J_k}(x)|.
\end{align*}
Note that this also holds if $J_ke=0$ for some $k$. Finally, squaring left- and right-hand side of the latter inequality and applying Cauchy-Schwarz, the claim follows from the fact that there is some $c>0$ such that $\gamma_k(x,H_x(x))\ge c$ for $x\in K$ and all $k=1,\ldots,N$.
\end{proof}
\end{document}